\documentstyle[amscd,amssymb,verbatim]{amsart}
\newcommand{\cohproj}{\operatorname{cohproj}}

\newcommand{\tr}{\operatorname{tr}}
\newcommand{\Nm}{\operatorname{Nm}}

\newcommand{\coker}{\operatorname{coker}}

\newcommand{\SL}{\operatorname{SL}}

\newcommand{\Coh}{\operatorname{Coh}}
\newcommand{\GG}{{\cal G}}
\newcommand{\CC}{{\cal C}}

\newcommand{\Aut}{\operatorname{Aut}}

\renewcommand{\Sp}{\operatorname{Sp}}

\newcommand{\Proj}{\operatorname{Proj}}

\newcommand{\si}{\sigma}

\newcommand{\Pic}{\operatorname{Pic}}

\newcommand{\Id}{\operatorname{Id}}
\newcommand{\ga}{\gamma}
\newcommand{\de}{\delta}

\renewcommand{\ker}{\operatorname{ker}}

\numberwithin{equation}{section}

\newtheorem{thm}{Theorem}[section]
\newtheorem{prop}[thm]{Proposition}
\newtheorem{lem}[thm]{Lemma}
\newtheorem{cor}[thm]{Corollary}

\newenvironment{rem}{\vspace{3mm}\noindent
{\bf Remark.}}{\vspace{3mm}}
\newenvironment{rems}{\vspace{3mm}
\noindent {\bf Remarks.}}{\vspace{3mm}}
\newenvironment{ex}{\vspace{3mm}\noindent
{\bf Example.}}{\vspace{3mm}}
\newenvironment{exs}{\vspace{3mm}\noindent
{\bf Examples.}}{\vspace{3mm}}

\newcommand{\Pf}{\noindent {\it Proof}}
\newcommand{\id}{\operatorname{id}}

\newcommand{\ov}{\overline}

\renewcommand{\Im}{\operatorname{Im}}

\newcommand{\rk}{\operatorname{rk}}
\newcommand{\ra}{\rightarrow}

\newcommand{\FF}{{\cal F}}

\newcommand{\PP}{{\cal P}}

\newcommand{\SS}{{\cal S}}

\newcommand{\OO}{{\cal O}}

\newcommand{\Hom}{\operatorname{Hom}}

\newcommand{\End}{\operatorname{End}}

\renewcommand{\a}{\alpha}
\renewcommand{\b}{\beta}

\newcommand{\De}{\Delta}

\newcommand{\th}{\theta}
\newcommand{\C}{{\Bbb C}}
\newcommand{\R}{{\Bbb R}}
\newcommand{\Z}{{\Bbb Z}}
\newcommand{\Q}{{\Bbb Q}}

\newcommand{\Ga}{\Gamma}

\newcommand{\sub}{\subset}
\newcommand{\ed}{\qed\vspace{3mm}}

\title{Noncommutative two-tori with real multiplication
as noncommutative projective varieties}
\author{A.~Polishchuk}
\address{Department of Mathematics and Statistics,
Boston University, Boston, MA 02215}
\email{apolish@@math.bu.edu}
\thanks{This work was partially supported by NSF grant DMS-0070967}

\begin{document}
\begin{abstract}
We define analogues of homogeneous coordinate algebras for
noncommutative two-tori with real multiplication.
We prove that the categories of
standard holomorphic vector bundles on such noncommutative
tori can be described
in terms of graded modules over appropriate
homogeneous coordinate algebras.
We give a criterion for such an algebra to be Koszul and prove
that the Koszul dual algebra also comes from some
noncommutative two-torus with real multiplication.
These results are based on the techniques of \cite{PS}
allowing to interpret all the data in terms
of autoequivalences of the derived categories of coherent
sheaves on elliptic curves.
\end{abstract}

\maketitle

\centerline{\sc Introduction}

\bigskip

Noncommutative algebraic
geometry is usually understood as the study of certain abelian
categories replacing the usual category of (quasi-)coherent sheaves
(see \cite{AZ}, \cite{KR}, \cite{Ros}).
For example, noncommutative projective schemes correspond
to certain categories defined in terms of modules over graded algebras
in the way analogous to Serre's theorem (see \cite{AZ}).
However, it is rather disappointing that at present there is almost no
connection between noncommutative algebraic varieties over $\C$
and noncommutative topological spaces, which according
to Connes \cite{Connes} are described by $C^*$-algebras.
One of the indications that such a connection exists
is provided by the work \cite{CDV}, where Sklyanin algebras are
related to some noncommutative manifolds.
In the present paper we give another example of a relation of this kind.
Namely, we show that noncommutative two-tori admitting
``real multiplication" (i.e., nontrivial Morita autoequivalences) can be
considered as underlying noncommutative topological spaces for
certain noncommutative algebraic varieties.

This relation is not so surprizing given the recent studies of complex
geometry on noncommutative two-tori (see \cite{S},\cite{DS},\cite{PS}).
It may only seem a little odd that real multiplication
is relevant for our picture. Let us briefly explain this.
Recall that the homogeneous coordinate algebra of a projective scheme 
is defined using tensor powers of an ample line bundle. In noncommutative
world one can only take tensor powers of a bimodule. Therefore,
in order to construct an analogue of such algebra for a 
noncommutative two-torus $T_{\th}$, where $\th\in\R$,
one has to find a bimodule over the ring of functions on $T_{\th}$
which would be ample in approriate sense. The natural choice
would be one of the so called basic modules.
Now the Morita equivalence theory for noncommutative two-tori implies that
an interesting bimodule can be found among basic modules only when
the parameter $\th$ is stabilized by a nontrivial element
of $\SL_2(\Z)/\{\pm 1\}$ under the fractional-linear action of this
group on $\R\cup\{\infty\}$. In other words, the category of vector
bundles on $T_{\th}$ should have a nontrivial Morita autoequivalence.
Note however, that there exists a generalization of
the standard approach to noncommutative projective schemes in which
graded algebras are replaced by more general objects called $\Z$-algebras
(see \cite{BP}, \cite{SVdB}). If one allows these more general
noncommutative ``$\Z$-projective schemes" then the condition that
$T_{\th}$ has real multiplication becomes unnecessary.


The results of this paper depend heavily on
the study of categories of holomorphic vector bundles on
$T_{\th}$ in \cite{PS}. Recall that in {\it loc.~cit.} we considered
only certain class of
holomorphic bundles on $T_{\th}$ that we called {\it standard} and
we constructed a fully faithful functor
from the category of such bundles to
the derived category $D^b(X)$ of coherent sheaves on some
elliptic curve $X$. Moreover, we proved that the image of this functor
consists of stable objects in the heart $\CC^{\th}$ of certain 
nonstandard $t$-structure on $D^b(X)$ associated with $\th$
(see \ref{tstrsec}; these $t$-structures were defined in \cite{Bridge}).
We conjecture that every
holomorphic bundle on $T_{\th}$ is a successive extension of
standard holomorphic bundles. If true, this would imply an equivalence
of $\CC^{\th}$ with the category of all holomorphic bundles on $T_{\th}$
(for irrational $\th$). Since we do not know the validity of this
conjecture, we simply replace the category of all holomorphic bundles
on $T_{\th}$ by $\CC^{\th}$.
This allows us to switch from the context of noncommutative complex
geometry on 
$T_{\th}$ to the study of the $t$-structure on $D^b(X)$ associated with $\th$.
Nontrivial Morita autoequivalences of $T_{\th}$
appearing when $\th$ is a quadratic irrationality
correspond to nontrivial autoequivalences $F:D^b(X)\ra D^b(X)$ preserving
the corresponding $t$-structure. 

The graded algebras associated with $T_{\th}$ can now be viewed as
examples of the following general construction.
Given an additive category $\CC$, an additive
functor $F:\CC\ra\CC$ and an object $O$
of $\CC$, we define an associative graded ring
$$A_{F,O}=\oplus_{n\ge 0}\Hom_{\CC}(O,F^n(O)),$$
where $(F^n:\CC\ra\CC, n\ge 0)$
are the functors obtained by iterating $F$ (so $F^0=\Id_{\CC}$).
The multiplication is defined as the composition of
the natural maps
\begin{align*}
&\Hom_{\CC}(O,F^m(O))\otimes\Hom_{\CC}(O,F^n(O))
\ra\Hom_{\CC}(F^n(O),F^{m+n}(O))\otimes\Hom_{\CC}(O,F^n(O))\\
&\ra
\Hom_{\CC}(O,F^{m+n}(O)).
\end{align*}
The homogeneous coordinate ring of a projective scheme $X$ appears as a 
particular case of this construction 
when $\CC$ is the category of coherent sheaves on $X$, 
$F$ is the functor of tensoring with an ample
line bundle $L$ on $X$, $O=\OO_X$ is the structure sheaf.
Slightly more general rings are obtained when taking $F$
to be of the form $F(A)=L\otimes\si^*A$, where $\si$ is
an automorphism of $X$. The corresponding rings are twisted
homogeneous coordinate rings considered in \cite{AVdB}.

The example relevant for noncommutative tori with real multiplication
is when $\CC=\CC^{\th}\sub D^b(X)$, where $X$ is an elliptic curve,
$F$ is the autoequivalence of $D^b(X)$ preserving $\CC^{\th}$.
In section \ref{algsec} we study corresponding graded algebras $A_{F,\FF}$,
where $\FF$ is a stable object of $\CC^{\th}$.
Namely, we compute the Hilbert series of $A_{F,\FF}$ and formulate
simple criterions in terms of discrete invariants of $(F,\FF)$ for
the algebra
$A_{F,\FF}$ to be generated in degree $1$, to be quadratic, and to be Koszul.
We also observe that if $A_{F,\FF}$ is Koszul then the Koszul dual algebra
is again of the same form: it is isomorphic to $A_{R_{\FF}\circ F^{-1},\FF}$,
where $R_{\FF}$ is certain twist functor associated with $\FF$ 
(see section \ref{Koszulsec}). 

In section \ref{projsec} we prove that every category $\CC^{\th}$, where
$\th$ is a quadratic irrationality, contains an {\it ample} sequence
of objects of the form $(F^n \FF)$, where $F:\CC^{\th}\ra\CC^{\th}$ is
an autoequivalence. This means that $\CC^{\th}$ can be
recovered from the corresponding graded algebra $A_{F,\FF}$ by the
noncommutative analogue of $\Proj$-construction considered in \cite{AZ}.
One technical point is that the categories $\CC^{\th}$ are non-Noetherian,
so we have to apply the main result of \cite{P-coh} that generalizes (a part of) 
the main theorem of \cite{AZ} to non-Noetherian case.

It would be interesting to try to extend some of our results
to more general algebras of the form $A_{F,\FF}$, where $F$ is an 
autoequivalence of the derived category $D^b(X)$ of coherent sheaves on a smooth
projective variety $X$, $\FF$ is an object of $D^b(X)$. The first natural question
is whether there are interesting examples when $F$ preserves some $t$-structure
on $D^b(X)$. In the case when $X$ is an abelian variety one source of such examples
should be given by noncommutative tori generalizing the picture described in \cite{PS}.

Another perspective for the future work is
to try to connect our results with Manin's program in \cite{Manin} 
to use noncommutative two-tori with real multiplication for
the explicit construction of the maximal abelian extensions
of real quadratic fields.

\noindent
{\it Convention}. With the exception of section \ref{torisec} all the objects
(varieties, categories) are defined over an arbitrary field $k$.

\section{Preliminaries on
derived categories of elliptic curves}\label{autoeqsec}

\subsection{Structure of the group of autoequivalences}\label{groupsec}

Let $X$ be an elliptic curve,
$\Aut(D^b(X))$ be the group of (isomorphism classes of)
exact autoequivalences of $D^b(X)$.
There is a natural surjective homomorphism
$$\pi:\Aut(D^b(X))\ra\SL_2(\Z),$$
defined by the rule $\pi(F)=g\in\SL_2(\Z)$,
such that for every object $\FF\in D^b(X)$ one has
$$\left(\matrix\deg F(\FF) \\ \rk F(\FF)\endmatrix\right)=
g\left(\matrix \deg \FF \\ \rk \FF\endmatrix\right).$$
For example, if $F$ is the functor of tensoring with a line
bundle $L$ then $F$ projects to the matrix
$\left(\matrix 1 & \deg(L) \\ 0 & 1\endmatrix\right)$.
Let $\SS:D^b(X)\ra D^b(X)$ be the Fourier-Mukai transform
considered as an autoequivalence of $D^b(X)$ via the
isomorphism $\hat{X}\simeq X$. Then
$\pi(\SS)=\left(\matrix 0 & -1\\ 1 & 0\endmatrix\right)$.
The shift functor $\FF\mapsto \FF[1]$ maps under $\pi$ to the matrix
$-\id\in\SL_2(\Z)$. We denote by $\Aut(X)$
the group of automorphisms of $X$ preserving the neutral
element. It can be identified with a subgroup of $\Aut(D^b(X))$: to every
automorphism $\si:X\ra X$ there corresponds an autoequivalence
$\si_*=(\si^{-1})^*:D^b(X)\ra D^b(X)$. Clearly, the homomorphism
$\pi$ is trivial on this subgroup.

On the other hand, for every abelian variety $X$ there is a homomorphism
$$\ga_X:\Aut(D^b(X))\ra\Sp(X\times\hat{X}),$$
where $\Sp(X\times\hat{X})$ is the group of
{\it symplectic} automorphisms of $X\times\hat{X}$, i.e.
automorphisms preserving the line bundle $p_{14}^*\PP\otimes
p_{23}^*\PP^{-1}$ on $(X\times\hat{X})^2$, where $\PP$ is the Poincar\'e
line bundle on $X\times\hat{X}$.
The homomorphism $\ga_X$ was defined by
Orlov in \cite{O}, Cor. 2.16. He also proved that it fits into an exact
sequence
$$1\ra (X\times\hat{X})(k)\times\Z\ra 
\Aut(D^b(X))\stackrel{\ga_X}{\ra}\Sp(X\times\hat{X})\ra 1,$$
where the subgroup $(X\times\hat{X})(k)$
corresponds to functors of translation by points of $X$ and of
tensor products with line bundles in $\Pic^0(X)$, while
$\Z\sub\Aut(D^b(X))$ is the subgroup of shifts $A\mapsto A[n]$. 
More precisely, to a point $(x,\xi)\in (X\times\hat{X})(k)$
one associates the autoequivalence
$$\Phi_{(x,\xi)}:D^b(X)\ra D^b(X):\FF\mapsto t^*_{-x}(\FF)\otimes
\PP|_{X\times\xi},$$
where $t_{x'}:X\ra X$ denotes the translation by $x'\in X(k)$.
The subgroup $(X\times\hat{X})(k)\sub\Aut(D^b(X))$ is normal
and the adjoint action of $F\in\Aut(D^b(X))$ is given precisely
by $\ga_X(F)$ (see \cite{O}, Cor. 2.13).

In the case of an elliptic curve
we can identify $X$ with $\hat{X}$, so the group
$\Sp(X\times\hat{X})$ can be identified with the group
$\Sp(X\times X)$ of matrices
$\left(\matrix a & b\\ c & d\endmatrix \right)$ with entries
in the ring $\End(X)$ satisfying the equations 
$$\ov{a}d-\ov{c}b=1,\ a\ov{d}-b\ov{c}=1$$
$$\ov{a}c=\ov{c}a,\  \ov{b}d=\ov{d}b,\
a\ov{b}=b\ov{a},\ c\ov{d}=d\ov{c},
$$
where $f\mapsto\ov{f}$ is the Rosati involution on $\End(X)$.

\begin{lem} The group $\Sp(X\times X)$ is isomorphic
to $(\Aut(X)\times\SL_2(\Z))/\{\pm 1\}$, where the subgroup
$\{\pm 1\}$ is embedded diagonally.
\end{lem}

\Pf . Note that for every $a\in\Aut(X)$ one has
$\ov{a}a=1$. Hence, $\Aut(X)$ embeds into $\Sp(X\times X)$
as the central subgroup of diagonal matrices 
$\left(\matrix a & 0 \\ 0 & a\endmatrix\right)$.
We have to prove that every element in $\Sp(X\times X)$ is a product
of such a matrix with a matrix in $\SL_2(\Z)$. If one of the entries
of a matrix $\left(\matrix a & b\\ c & d\endmatrix \right)\in\Sp(X\times X)$
is zero then this is easy. Assuming that all the entries are non-zero
we note that the condition $a\ov{b}\in\Z$ implies
that $a\in\Q b$. Therefore, we can write
$a=a'r$, $b=b'r$ for some $r\in\End(X)_{\Q}$ and a pair of relatively
prime integers $(a',b')$. From the condition $a,b\in\End(X)$ we
immediately derive that $r\in\End(X)$. Using the conditions 
$\ov{a}c\in\Z$, $\ov{b}d\in\Z$ we can also write
$c=c'\ov{r}^{-1}$, $d=d'\ov{r}^{-1}$ for some rational numbers $(c',d')$.
Since $c$, $d$, and $\ov{r}$ are elements of $\End(X)$ we obtain
that $c'$ and $d'$ are integers. The equation
$\ov{a}d-\ov{c}b=1$ implies that the matrix 
$\left(\matrix a' & b'\\ c' & d'\endmatrix \right)$ belongs to
$\SL_2(\Z)$. Finally, since $c'\ov{r}^{-1}\in\End(X)$ and
$d'\ov{r}^{-1}\in\End(X)$ it follows that $\ov{r}^{-1}\in\End(X)$,
so $r$ is a unit.
\ed

Let us set $\ov{\Aut}(D^b(X))=\Aut(D^b(X))/(X\times\hat{X})(k)$.
According to the above lemma the homomorphism $\ga_X$ induces a
surjective homomorphism
$$\ga_X:\ov{\Aut}(D^b(X))\ra (\Aut X\times\SL_2(\Z))/\{\pm 1\}$$
with the kernel $\Z$.
The homomorphism $\pi$ also factors through a homomorphism
$$\pi:\ov{\Aut}(D^b(X))\ra\SL_2(\Z).$$
It is easy to check that the
homomorphisms $\ov{\Aut}(D^b(X))\ra\SL_2(\Z)/\{\pm 1\}$
induced by $\ga_X$ and $\pi$ differ by an automorphism of
$\SL_2(\Z)$. On the other hand,
there is a natural action of $\Aut(D^b(X))$ on $K_0(X)$
that preserves the subgroup $K_0(X)_0\sub K_0(X)$
consisting of classes of zero degree and zero rank.
Note that the determinant
gives an isomorphism $\det:K_0(X)_0\ra\Pic^0(X)$.
From this one can see that the action of $\Aut(D^b(X))$ on $K_0(X)_0$
factors through $\ov{\Aut}(D^b(X))$.

\begin{thm} There exists a homomorphism
$\ov{\Aut}(D^b(X))\ra\Aut(X):F\mapsto\a_F$, such that
$$F(a)=\a_{F}(a)$$
for every $a\in K_0(X)_0$.
The induced homomorphism
$$\ov{\Aut}(D^b(X))\ra\Aut(X)\times\SL_2(\Z):F\mapsto (\a_F,\pi(F))$$
fits into the exact sequence
$$1\ra 2\Z\ra \ov{\Aut}(D^b(X))\ra\Aut(X)\times\SL_2(\Z)\ra 1,$$
where $2\Z$ is the subgroup of even shifts $A\mapsto A[2n]$.
\end{thm}

\Pf . We have an exact sequence
\begin{equation}\label{autexseq}
1\ra\Z\ra\ov{\Aut}(D^b(X))/\Aut(X)\stackrel{\ov{\pi}}{\ra}
\SL_2(\Z)/\{\pm 1\}\ra 1,
\end{equation}
where the homomorphism $\ov{\pi}$ is induced by $\pi$, $\Z$
is the subgroup of shifts.
Indeed, this follows from the fact that $\ov{\pi}$ differs
from the homomorphism induced by $\ga_X$ by an automorphism
of $\SL_2(\Z)/\{\pm 1\}$ and from the fact that $\ker(\ga_X)=\Z$.
One consequence of this is that
$\ov{\Aut}(D^b(X))$ is generated by the subgroup $\Aut(X)$
together with the shift $A\mapsto A[1]$,
the tensoring functor $T(A)=A\otimes L$,
where $L$ is a line bundle of degree $1$,
and the Fourier-Mukai transform $\SS$. Since the action of all of these
autoequivalences on $K_0(X)$ is known, we derive that for every
$F\in\ov{\Aut}(D^b(X))$ there exists an automorphism $\a_F\in\Aut(X)$
such that $F$ acts on $K_0(X)_0$ in the same way as $\a_F$.
Similarly, we can consider the action of $F$ on
$K_0(X\otimes_k\ov{k})_0$, where $\ov{k}$ is an algebraic closure of
$k$. The above argument shows
that there exists unique $\a_F$ defined over
$k$, such that $F$ acts on $K_0(X\otimes_k\ov{k})_0$ as $\a_F$.
It is clear that the homomorphism $F\mapsto\a_F$ restricts to the
identity map on $\Aut(X)\sub\ov{\Aut}(D^b(X))$. Together with surjectivity
of $\pi$ this immediately implies surjectivity of the
map $F\mapsto (\a_F,\pi(F))$.
The kernel of this map clearly contains the subgroup of even shifts
$2\Z\sub\Z\sub\ov{\Aut}(D^b(X))$.
Using the exact sequence (\ref{autexseq}) one can
easily see that this kernel coincides with $2\Z$.
\ed

Let $\ov{\Aut}(D^b(X))^0\sub\ov{\Aut}(D^b(X))$
(resp., $\Aut(D^b(X))^0\sub\Aut(D^b(X))$) be the subgroup
consisting of $F$ with $\a_F=1$. From the above theorem
we get an exact sequence
$$1\ra 2\Z\ra\ov{\Aut}(D^b(X))^0\ra\SL_2(\Z)\ra 1.$$
The following proposition can be viewed as an analogue
of the theorem of the cube for autoequivalences of $D^b(X)$.

\begin{prop}\label{cubefunprop}
For every $F\in\Aut(D^b(X))$ one has
\begin{equation}\label{cubefun}
[F^3(a)]-(N_F+\a_F)[F^2(a)]+(1+N_F\a_F)[F(a)]-\a_F[a]=0
\end{equation}
in $K_0(X)$ for every $a\in K_0(X)$, where
$N_F=\tr(\pi(F))$.
In particular, if $F\in\Aut(D^b(X))^0$ then
$$[F^3(a)]-(N_F+1)[F^2(a)]+(N_F+1)[F(a)]-[a]=0.$$
\end{prop}

\Pf . Set $N=N_F$.
We claim that $[F^2(a)]-N[F(a)]+[a]\in K_0(X)_0\sub K_0(X)$
for every $a\in K_0(X)$. Indeed, this follows immediately from
the fact that $F$ acts as an element $g$ on the pair $(\deg,\rk)$,
where $g^2-Ng+1=0$. It remains to apply $F$ to this element
in $K_0(X)_0$ and to use the definition of $\a_F$.
\ed

\begin{ex} Let $F$ be of the form $(\otimes L)\circ\si^*$,
where $L$ is a line bundle, $\si$ is a translation by a point on $X$.
Then $\a_F=1$ and
$$F^2\simeq (\otimes L\otimes\si^*L)\circ\si^*,$$
$$F^3\simeq (\otimes L\otimes\si^*L\otimes(\si^2)^*L)\circ\si^*,$$
so (\ref{cubefun}) amounts to the identity
$$[L\otimes \si^*L\otimes(\si^2)^*L]-3[L\otimes\si^*L]
+3[L]-[\OO_X]=0,$$
or equivalently
$$(\si^2)^*L\simeq \si^*L^2\otimes L^{-1}.$$
On the other hand, this isomorphism is a direct consequence
of the theorem of the cube. This is why we can view Proposition
\ref{cubefunprop} as an analogue of this theorem for more general
autoequivalences.
\end{ex}

\subsection{$t$-structures on $D^b(X)$}\label{tstrsec}

Let us say that an object $\FF\in D^b(X)$ is {\it stable}
if $\FF=V[n]$, where $n\in\Z$, $V$ is either a stable
vector bundle or the structure sheaf of a $k$-point.
It is easy to see that an object $\FF\in D^b(X)$ is stable if and only if
$\Hom(\FF,\FF)=k$ (since every object in $D^b(X)$ is isomorphic to the direct
sum of its cohomology sheaves). 

Below we are going to use some basic notions and results
of the torsion theory that can be found in \cite{HRS}.
For every real number $\th$ we consider a $t$-structure
$(D^{\th,\le 0},D^{\th,\ge 0})$ on $D^b(X)$ defined as follows. 
First, let us define a torsion pair $(\Coh_{>\th},\Coh_{\le\th})$ in
the category $\Coh(X)$ of coherent sheaves on $X$.
By the definition, $F\in\Coh_{>\th}$ (resp. $\Coh_{\le\th}$)
if all semistable factors of $F$ have slope $>\th$ (resp. $\le\th$).
Note that objects of $\Coh_{>\th}$ are allowed to have arbitrary torsion
(we consider torsion sheaves as having the slope $+\infty$).
Now the $t$-structure associated with $\th$ is defined by the rule
$$D^{\th,\le 0}:=\{ K\in D^b(X):\ H^{>0}(K)=0,\ H^0(K)\in\Coh_{>\th}\},$$
$$D^{\th,\ge 1}:=\{ K\in D^b(X):\ H^{<0}(K)=0,\ H^0(K)\in\Coh_{\le\th}\}.$$ 
The fact that this is indeed a $t$-structure follows from the torsion
theory (see \cite{HRS}). The heart of this $t$-structure is
$\CC^{\th}(X):=D^{\th,\le 0}\cap D^{\th,\ge 0}$. It
is equipped with the torsion pair $(\Coh_{\le\th}[1],\Coh_{>\th})$.

It is convenient to
extend the above definition to $\th=\infty$ by letting
$(D^{\infty,\le 0},D^{\infty,\ge 0})$ to be the standard $t$-structure
on $D^b(X)$. 
In the following proposition we list some properties of these
$t$-structures.
Let us denote $v_{\FF}=(\deg(\FF),\rk(\FF))\in\Z^2$ for $\FF\in D^b(X)$.
Note that for $F\in\Aut(D^b(X))$ one has
$v_{F(\FF)}=\pi(F)(v_{\FF})$, where
$\SL_2(\Z)$ acts on $\Z^2$ as on column vectors.

\begin{prop}\label{cohdim}
(i) The category $\CC^{\th}(X)$ has cohomological dimension
$1$ and there is an equivalence $D^b(\CC^{\th}(X))\simeq D^b(X)$.

\noindent
(ii) If $\th\in\R\setminus\Q$ then one has
$$\{v_{\FF} |\ \FF\in\CC^{\th}(X),\FF\neq 0\}=H_{\th}\cap\Z^2,$$
where $H_{\th}=\{(x_1,x_2)|\ x_1-\th x_2>0\}$.
If $\th\in\Q$ then the above statement is
true with the half-plane $H_{\th}$ replaced by its union
with the ray $\R_{\le 0}(\th,1)$. 
\end{prop}

\Pf . (i) Let us first prove that $\Hom_{D^b(X)}^{>1}(A,B)=0$ for every
$A,B\in\CC^{\th}(X)$. If both objects
$A$ and $B$ belong to one of the subcategories
$\Coh_{>\th}$ or $\Coh_{\le \th}[1]$, then the assertion is clear.
If $A\in\Coh_{>\th}$, $B\in\Coh_{\le\th}[1]$ then
$\Hom_{D^b(X)}^i(A,B)=\Hom_{D^b(X)}^{i+1}(A,B[-1])=0$ for
$i\ge 1$. On the other hand,
$$\Hom_{D^b(X)}^i(B,A)=\Hom_{D^b(X)}^{i-1}(B[-1],A)\simeq
\Hom_{D^b(X)}^{2-i}(A,B[-1])^*=0$$
for $i\ge 2$, since $\Hom_{\Coh(X)}(A,B[-1])=0$.

The second assertion follows from this by the standard argument
(see e.g. \cite{Beil}).

\noindent
(ii) This follows from the fact that primitive lattice vectors
contained in $H\cup\R_{\le 0}(\th,1)$ are exactly vectors 
$v_{\FF}$, where $\FF$ is a stable object belonging to
$\Coh_{>\th}$ or $\Coh_{\le\th}[1]$.
\ed

It is not difficult to calculate the action of autoequivalences of $D^b(X)$
on these $t$-structures. We denote by $\th\mapsto g\th=
\frac{a\th+b}{c\th+d}$, where
$g=\left(\matrix a & b \\ c & d\endmatrix\right)\in\SL_2(\Z)$, the standard
fractional-linear action of $g$ on $\R\cup\{\infty\}$.

\begin{prop}\label{tstrprop}
For every $F\in\Aut(D^b(X))$ and every $\th\in\R\cup\{\infty\}$
one has $F(D^{\th,\le 0})=D^{g\th,\le 0}[n]$
(resp., $F(D^{\th,\ge 0})=D^{g\th,\ge 0}[n]$),
where $g=\pi(F)\in\SL_2(\Z)$, $n$ is some integer.
\end{prop}

\Pf . The assertion is clear when $F$ is a translation or the tensor product
with a line bundle, or the pull-back under an automorphism. Hence,
it suffices to consider the case $F=\SS$, where $\SS$ is the Fourier-Mukai
transform. For every segment $I\sub\R\cup\{+\infty\}$ let us denote by
$\Coh_I(X)\sub\Coh(X)$ the full subcategory in $\Coh(X)$ consisting
of sheaves $F$ such that all semistable factors of $F$ have slope in $I$
(recall that torsion sheaves have slope $+\infty$). 
Assume first that $\th\in\R$, $\th>0$. Since $\SS$ transforms the slopes
by the map $\mu\mapsto-\mu^{-1}$, we have
$$\SS(\Coh_{(\th,+\infty]})=\Coh_{(-\th^{-1},0]},$$
$$\SS(\Coh_{(0,\th]}[1])=\Coh_{(-\infty,-\th^{-1}]}[1],$$
$$\SS(\Coh_{(-\infty,0]}[1])=\Coh_{(0,+\infty]}.$$
This immediately implies that $\SS$ sends the $t$-structure associated
with $\th$ to the $t$-structure associated with $-\th^{-1}$.
Since $\SS^2=[-\id]^*[-1]$, it follows that for $\th\in\R$, $\th<0$ 
one has $\SS(D^{\th,\le 0},D^{\th,\ge 0})=
(D^{-\th^{-1},\le 0}[-1],D^{-\th^{-1},\ge 0}[-1])$.
Similarly, it easy to see that $\SS$ switches (up to a shift)
the standard $t$-structure with the $t$-structure corresponding to $\th=0$.
\ed

Let us consider the bilinear form $\chi(\FF_1,\FF_2):=\sum_i (-1)^i\dim\Hom^i(\FF_1,\FF_2)$
on $K_0(X)$, where 
$\Hom^i(\FF_1,\FF_2):=\Hom_{D^b(X)}(\FF_1,\FF_2[i])$. It is easy to see
that 
$$\chi(\FF_1,\FF_2)=-\det(v_{\FF_1},v_{\FF_2}).$$
The kernel of $\chi$ is exactly the subgroup $K_0(X)_0\sub K_0(X)$
consisting of elements of zero degree and zero rank.
Abusing the notation we also set $\chi(v,v'):=-\det(v,v')$ for $v,v'\in\Z^2$,
so that $\chi(\FF_1,\FF_2)=\chi(v_{\FF_1},v_{\FF_2})$.

The following lemma generalizes to categories $\CC^{\th}(X)$ the well-known fact
about stable bundles on $X$.

\begin{lem}\label{nzhomlem} Let $\FF_1,\FF_2$ be a pair of 
stable objects in $\CC^{\th}(X)$ such that
$\chi(\FF_1,\FF_2)>0$. Then $\Hom^1(\FF_1,\FF_2)=0$ and
$\dim\Hom(\FF_1,\FF_2)=\chi(\FF_1,\FF_2)$.
\end{lem}

\Pf . By Proposition \ref{cohdim}(i) we have
$\chi(\FF_1,\FF_2)=\dim\Hom(\FF_1,\FF_2)-\dim\Hom^1(\FF_1,\FF_2)$,
so it is enough to check the vanishing of $\Hom^1(\FF_1,\FF_2)$.
In the case when both $\FF_1$ and $\FF_2$ belong to $\Coh_{>\th}$
(resp., $\Coh_{\le\th}[1]$) the assumption $\chi(v_{\FF_1},v_{\FF_2})>0$
implies that
$\mu(\FF_1)<\mu(\FF_2)$ (resp., $\mu(\FF_1[-1])<\mu(\FF_2[-1])$).
Hence, in this case the assertion is clear. It is easy to see that the case 
$\FF_1\in\Coh_{\le\th}[1]$, $\FF_2\in\Coh_{>\th}$ cannot occur. 
Indeed, since the vectors $v_{\FF_1}$ and $v_{\FF_2}$ belong to $H_{\th}$
and $\chi(v_{\FF_1},v_{\FF_2})>0$, the condition $\rk\FF_2>0$ implies
that $\rk\FF_1>0$.
In the
remaining case $\FF_1\in\Coh_{>\th}$ and $\FF_2\in\Coh_{\le\th}[1]$, so 
the assertion follows from the vanishing of 
$\Hom^{2}(\FF_1,\FF_2[-1])$.
\ed

\section{Noncommutative tori, autoequivalences of $D^b(X)$ and
related algebras}\label{algsec}

\subsection{Morita autoequivalences of noncommutative tori
and analogues of homogeneous coordinate rings}\label{torisec}

We refer to \cite{Rief} for an introduction and a survey of main
results in the theory of noncommutative tori.
Recall that for every $\th\in\R$ the algebra 
$A_{\th}$ of smooth functions on the
noncommutative torus $T_{\th}$ is the algebra of
series $\sum a_{n_1,n_2} U_1^{n_1}U_2^{n_2}$ in variables
$U_1$, $U_2$ satisfying the relation
$$U_1U_2=\exp(2\pi i \th)U_2U_1,$$ 
such that the coefficient function
$(n_1,n_2)\ra a_{n_1,n_2}\in\C$ is rapidly decreasing at infinity.
By the definition, vector bundles on $T_{\th}$ are finitely generated
projective $A_{\th}$-modules (we always consider right modules).
A complex structure on $T_{\th}$ is given by a derivation
$\de_{\tau}:A_{\th}\ra A_{\th}$, such that $\de_{\tau}(U_1)=\tau$,
$\de_{\tau}(U_2)=1$, where $\tau\in\C\setminus\R$ (following \cite{PS}
we will usually impose the condition $\Im(\tau)<0$). We denote
by $T_{\th,\tau}$ the noncommutative torus $T_{\th}$ equipped with
this complex structure.
A holomorphic structure on a vector bundle is given by an operator
$\nabla:E\ra E$ on the corresponding projective right $A_{\th}$-module
satisfying $\nabla(ea)=\nabla(e)a+e\de_{\tau}(a)$, where $e\in E$,
$a\in A_{\th}$.
As in \cite{PS} we only consider {\it standard} holomorphic vector
bundles on $T_{\th,\tau}$ that are given by certain family of standard
holomorphic structures on {\it basic modules}.

Recall that a basic $A_{\th}$-module $E$ is uniquely determined by its rank which
is a primitive positive element in $\Z+\Z\th$ 
(we assume that $\th$ is irrational).
The algebra of endomorphisms of $E$ can be identified
with $A_{\th'}$ for some $\th'\in\R$ and the functor
$E'\mapsto E'\otimes_{A_{\th'}}E$ is a Morita equivalence 
between the categories of right $A_{\th'}$-modules and right
$A_{\th}$-modules.
It is known that $\th'$ is necessarily of the form $\th'=g\th$
for some $g\in\SL_2(\Z)/\{\pm 1\}$. It is convenient to lift
elements of $\SL_2(\Z)/\{\pm 1\}$ to elements $g$ of $\SL_2(\Z)$
satsisfying the condition $c\th+d>0$, where
$g=\left(\matrix a & b \\ c & d\endmatrix\right)$. 
For every such $g\in\SL_2(\Z)$
there is a basic right $A_{\th}$-module $E_g(\th)$
of rank $\rk(E_g(\th))=c\th+d$
equipped with an isomorphism
$\End_{A_{\th}}(E_g(\th))\simeq A_{g\th}$ (see \cite{PS}, section 1.1).

Standard holomorphic structures on $E_g(\th)$ are
parametrized by the points of a complex elliptic curve 
$X_{\tau}=\C/(\Z+\Z\tau)$.
As we have shown in \cite{PS} (Prop. 3.1), fixing such a structure we obtain
the Morita equivalence $E'\mapsto E'\otimes_{A_{\th}}E_g(\th)$
between the categories of holomorphic bundles on $T_{g\th,\tau}$ and
$T_{\th,\tau}$ (preserving the subcategories of standard holomorphic
bundles).

We are interested in the situation when $g(\th)=\th$
for some nontrivial $g\in\SL_2(\Z)/\{\pm 1\}$.
It is easy to see that this happens exactly
when either $\th$ is rational or $\th$
generates a real quadratic extension of $\Q$.
In this case 
$E=E_g(\th)$ has a natural structure of $A_{\th}$-$A_{\th}$-bimodule, so
we can consider its tensor powers
$$E^{\otimes n}:=E\otimes_{A_{\th}}\ldots\otimes_{A_{\th}}
E\ (n\ \text{times})\simeq E_{g^n}(\th),$$
where the last isomorphism follows from 
the general formula for the tensor product of basic modules
(see e.g. Proposition 1.2 of \cite{PS}).
If we equip $E$ with a standard holomorphic structure then
$E^{\otimes n}$ acquires the induced holomorphic structure, so we can
consider the corresponding space of holomorphic vectors $H^0(E^{\otimes n})$.
Note that in order for these spaces to be non-zero we need to impose the
condition $c>0$ (where we assume that $\Im(\tau)<0$;
see \cite{PS}, Prop. 2.5).
For $n=0$ we set $E^{\otimes 0}=A_{\th}$ and equip it
with the standard holomorphic structure $\de_{\tau}$. 
Now there is a natural structure of an associative algebra on
$$B_E:=\oplus_{n\ge 0}H^0(E^{\otimes n})$$
given by the tensor product of holomorphic vectors.
Clearly, this algebra is a direct generalization of the homogeneous
coordinate ring. One can also  define analogues of twisted homogeneous
coordinate ring by changing the bimodule structure on $E$. Namely,
one can leave the right $A_{\th}$-module structure the same and
twist the left $A_{\th}$-module structure by some holomorphic
(i.e., commuting with $\de_{\tau}$) automorphism of $A_{\th}$.

Recall that in \cite{PS} we constructed an equivalence between
the derived category
of standard holomorphic vector bundles on $T_{\th,\tau}$
and the full subcategory
of stable objects in the derived category 
$D^b(X)$ of coherent sheaves on the elliptic
curve $X=X_{\tau}=\C/(\Z+\tau\Z)$. This equivalence sends each
standard holomorphic bundle of rank $m\th+n$ to a stable object in $D^b(X)$
of degree $m$ and rank $n$. The image of the category of holomorphic bundles
under this equivalence belongs to $\CC^{-\th^{-1}}$ (up to a shift).
Moreover, the Morita autoequivalence
$E'\mapsto E'\otimes_{A_{\th}}E_g(\th)$, where $g(\th)=\th$,
corresponds to some autoequivalence $F:D^b(X)\ra D^b(X)$ 
preserving $\CC^{-\th^{-1}}\sub D^b(X)$, such that 
$\pi(F)=g^t$ (the transposed matrix to $g$)
Note that this is compatible with Proposition \ref{tstrprop}
since $g^t(-\th^{-1})=-\th^{-1}$. It is also easy to see from the explicit
formulas for the equivalence of \cite{PS} that $F$ belongs
to the subgroup $\Aut(D^b(X))^0\sub\Aut(D^b(X))$ introduced
in \ref{groupsec}. By twisting the left
$A_{\th}$-module structure on $E_g(\th)$
with some holomorphic automorphisms
of $A_{\th}$ we can get more general autoequivalences $F$ with $\pi(F)=g^t$.

Let $(\FF_n, n\ge 0)$ be the image of the sequence of holomorphic bundles
$(E^{\otimes n})$ under the above equivalence of categories.
Then we have $\FF_n=F^n(\FF_0)$ and
there is a natural isomorphism of algebras
$$B_E\simeq A_{F,\FF_0}.$$
This isomorphism allows us to switch to the language of $t$-structures and
autoequivalences of $D^b(X)$ in the further study of these algebras.

It is sometimes convenient to change the point of view slightly. Namely,
the condition $g\th=\th$ 
is equivalent to the condition 
$$r^2-(a+d)r+1=0,$$
where $r=c\th+d$. In other words, $r$ is an eigenvalue of $g$.
Note that $\th$ can be recovered from the pair $(g,r)$ by
the formula $\th=(r-d)/c$.
Thus, we can start with an arbitrary matrix $g$ in $\SL_2(\Z)$ having
real positive eigenvalues and $c>0$. Then fixing
one of the eigenvalues $r$ of $g$, we get a family of graded algebras
$$B_{g,r}(\tau)=\oplus_{n\ge 0} H^0(T_{\th,\tau},E^{\otimes n})$$
parametrized by $\tau\in\C$ such that $\Im(\tau)<0$, 
where $\th=(r-d)/c$, $E=E_g(\th)$ is equipped with a standard holomorphic
structure $\ov{\nabla}_0$ (see \cite{PS}).
Note that for every upper-triangular matrix $u\in\SL_2(\Z)$ one has
$B_{ugu^{-1},r}(\tau)\simeq B_{g,r}(\tau)$,
so this family of algebras really 
lives on a double covering of $\SL_2(\Z)/Ad(\Z)$,
where $\Z$ is embedded as upper-triangular matrices in $\SL_2(\Z)$.

\begin{prop} There is an isomorphism of graded algebras
$$B_{g,r}(\tau)^{opp}\simeq B_{g',r^{-1}}(\tau),$$
where $g'=\left(\matrix d & b \\ c & a\endmatrix\right)$.
\end{prop}

\Pf . We have a natural isomorphism $A_{-\th}\simeq A_{\th}^{opp}$,
identical on generators $U_1$ and $U_2$. Thus, we can consider
$E=E_g(\th)$ as an $A_{-\th}$-$A_{-\th}$-bimodule. It is easy
to see that as such a bimodule $E$ is isomorphic to $E_{g'}(-\th)$.
Moreover, this isomorphism is compatible with standard holomorphic
structures up to a scalar. This implies the result.
\ed

\subsection{Algebras associated with autoequivalences and $t$-structures}


We want to look at the algebras $A_{F,\FF}$, where $\FF\in D^b(X)$
is a stable object, $F:D^b(X)\ra D^b(X)$ is an autoequivalence.
In order to get interesting algebras we would like to impose the
condition $F^n(\FF)\not\simeq\FF$ and
$\Hom(\FF,F^n(\FF))\neq 0$ for all
sufficiently large $n\ge 0$. We are going to show that
this condition implies that the corresponding element $\pi(F)\in\SL_2(\Z)$
has positive real eigenvalues.

\begin{lem}\label{eigenlem}
For $g\in\SL_2(\Z)$ and $v\in\Z^2\setminus\{0\}$
the following conditions are equivalent:

\noindent
(i) $\chi(v,g^n v)>0$ for all sufficiently large $n$;

\noindent
(i)' $\chi(v,g^n v)>0$ for all $n>0$;

\noindent
(ii) $g$ has positive real eigenvalues and $\chi(v,g v)>0$.
\end{lem}

\Pf . First, let us show that (i) implies (ii).
Set $N=\tr(g)$. Then $g^2-Ng+1=0$. Hence,
$$\chi(v,g^n v)-N\chi(v,g^{n-1} v)+\chi(v,g^{n-2} v)=0.$$
From this we immediately derive that
$$\sum_{n\ge 0} \chi(v,g^n v)t^n=\frac{Mt}{1-Nt+t^2},$$
where $M=\chi(v, g v)$.
Therefore, condition (ii) implies
that all coefficients of
the series $M(1-Nt+t^2)^{-1}$ except for a finite number
are positive. It is easy to see that for $N\ge 2$
(i.e., when both roots of the equation $t^2-Nt+1=0$ are real and positive)
all coefficients of the series $(1-Nt+t^2)^{-1}$ are positive. By
the change of variables $t\mapsto -t$ this implies that the series
$(1-Nt+t^2)^{-1}$ is alternating for $N\le -2$.
It is also easy to see directly that for $N=\pm 1$ or $N=0$ this series
still has infinitely many negative coefficients.
Hence, (i) implies that $M>0$ and
$N\ge 2$. The same argument shows that (ii) implies (i)'. 
\ed

\begin{prop}\label{eigenprop} 
Let $F:D^b(X)\ra D^b(X)$ be an autoequivalence
with $\pi(F)=g\in\SL_2(\Z)$, $\FF$ be a stable object of $D^b(X)$.
Then the following conditions are equivalent:

\noindent
(i) $F^n(\FF)\not\simeq\FF$ for $n>0$ and
$\Hom(\FF,F^n(\FF))\neq 0$ for all sufficiently large $n$;

\noindent
(ii) $g$ has positive real eigenvalues, $M=\chi(\FF, F(\FF))>0$,
and there exists
$\th\in\R\cup\{\infty\}$ such that $F$ preserves $\CC^{\th}$.

Under these conditions the Hilbert series of $A_{F,\FF}$ is equal to
$$H_{A_{F,\OO_X}}(t)=1+\frac{Mt}{1-Nt+t^2},$$
where $N=\tr(g)$.
Also, if $(\th,1)\in\R^2$ is an eigenvector of $g$
then $F$ preserves the subcategory $\CC^{\th}\sub D^b(X)$
(if $(1,0)$ is an eigenvector then we set $\th=\infty$).
\end{prop}

\Pf . (i)$\implies$ (ii). Note
that for every pair of non-isomorphic stable objects
$\GG,\GG'\in D^b(X)$ the graded space $\Hom^{\bullet}(\GG,\GG')$ is always
concentrated in one degree. Hence, (i) implies that
$\chi(\FF,F^n(\FF))=\chi(v_0,g^n v_0)>0$ for all 
sufficiently large $n$, where $v_0=v_{\FF}\in\Z^2$.
By Lemma \ref{eigenlem} this implies
that $g$ has positive eigenvalues and that $M>0$.
Also, from the proof of this lemma we obtain the formula
for the Hilbert series of $A_{F,\FF}$.
Now let $u\in\R^2$ be an eigenvector of $g$. Rescaling $u$
we can assume that either $u=(\th,1)$ or $u=(1,0)$.
In the latter case we set $\th=\infty$.
In either case we have $g(\th)=\th$ under the fractional-linear action.
Therefore, by Proposition \ref{tstrprop} we have
$F(\CC^{\th})=\CC^{\th}[m]$ for some $m\in\Z$.
Since $\FF$ belongs to $\CC^{\th}[i]$ for some $i\in\Z$, the nonvanishing
of $\Hom(\FF,F^n(\FF))$ for $n>>0$ implies that $m=0$.

(ii)$\implies$ (i). Let us consider vectors 
$$v_n=v_{F^n(\FF)}=g^n v_0\in\Z^2\sub\R^2.$$
Since $\chi(v_0,v_1)=M>0$, it follows
that $\chi(v_n,v_{n+1})>0$ for all $n\ge 0$.
Let $m$ be an integer such that $\FF\in\CC^{\th}[m]$.
Then all the vectors $v_n$ belong to the closed half-plane
$(-1)^m\ov{H_{\th}}\sub\R^2$ (see Proposition \ref{cohdim}).
Hence, $\chi(v_i,v_j)\ge 0$ for $i<j$. Moreover,
for a pair of stable objects $\GG,\GG'\in\CC^{\th}$ the vectors
$v_{\GG}$ and $v_{\GG'}$ can be proportional only if
$v_{\GG}=v_{\GG'}$. Hence, we have
$\chi(v_i,v_j)>0$ for $i<j$.
It remains to apply Lemma \ref{nzhomlem}. 
\ed


\begin{rems} 1. It is easy to deduce from the proof 
that for a pair $(F,\FF)$ such that $g=\pi(F)\in\SL_2(\Z)$
has positive real eigenvalues and $M>0$, the conditions
of the above proposition will be satisfied for $(F[n],\FF)$
for some $n\in 2\Z$.

\noindent 2. If $F$ satisfies the equivalent
conditions of the above proposition,
$g$ has two distinct eigenvalues and $(\th_1,1)$, $(\th_2,1)$
are the corresponding eigenvectors, then $F$ preserves both
subcategories $\CC^{\th_1}$ and $\CC^{\th_2}$. Hence, $F$ also preserves
$\CC^{\th_1}\cap\CC^{\th_2}$. Moreover, it is easy to see that
$\CC^{\th_1}\cap\CC^{\th_2}$ is a ``half'' of the natural $F$-invariant
torsion theory in each of the categories $\CC^{\th_1}$, $\CC^{\th_2}$
and that these categories are tiltings of each other with respect to these
torsion theories.
\end{rems}

One can rewrite the Hilbert series of $A_{F,\FF}$ as follows:
$$H_{A_{F,\FF}}(t)=\frac{1+(M-N)t+t^2}{1-Nt+t^2}.$$
In particular, we notice that the series
$H_{A_{F,\FF}}(-t)^{-1}$ has similar form but with $N$ and $M-N$
switched. Recall that if a graded algebra $A$ is Koszul then one has
$H_{A^!}(t)=H_A(-t)^{-1}$, where $A^!$ is the Koszul dual algebra.
Below we will show that under appropriate
conditions the algebra $A_{F,\FF}$ is Koszul with the dual also of
the form $A_{F',\FF}$.

\subsection{Koszul duality}\label{Koszulsec}

For every stable object $\FF\in D^b(X)$ we denote by $R_{\FF}$ the
the right twist corresponding to $\FF$. This is an
autoequivalence on $D^b(X)$ such that one has exact triangles
$$\GG\ra\Hom^{\bullet}(\GG,\FF)^*\otimes \FF\ra R_{\FF} \GG\ra \GG[1]$$
for all $\GG\in D^b(X)$ (see \cite{ST}; our notation is slightly
different). The quasi-inverse
autoequivalence is the left twist $L_{\FF}$, such that one
has an exact triangle
$$L_{\FF}(\GG)\ra\Hom^{\bullet}(\FF,\GG)\otimes\FF\ra\GG\ra L_{\FF}(\GG)[1].$$

\begin{thm}\label{Koszulthm}
Let $(F,\FF)$ be a pair 
satisfying the equivalent conditions of Proposition \ref{eigenprop}.
Let $\pi(F)=g\in\SL_2(\Z)$, $N=\tr(g)$,
$M=\chi(\FF,F(\FF))$,
and let $A=A_{F,\FF}$ be the corresponding graded algebra.
Then

\noindent
(a) $A$ is generated by $A_1$ over $A_0=k$ if and only if
$M\ge N+1$, or $M=N$ and
$$\det F^2(\FF)\not\simeq(\det F(\FF))^N\otimes\det(\FF)^{-1};$$

\noindent
(b) $A$ is a quadratic algebra if and only if $M\ge N+2$, or
$M=N+1$ and
$$\det F^3(\FF)\not\simeq(\det F^2(\FF))^{N+1}\otimes
(\det F(\FF))^{-N-1}\otimes\det(\FF);
$$

\noindent
(c) $A$ is Koszul if and only if $M\ge N+2$.
Moreover, in this case one has the following isomorphism
for the quadratic dual algebra:
$$A^!\simeq A_{R_{\FF}\circ F^{-1},\FF}.$$
\end{thm}

Let $(r,r^{-1})$ be the eigenvalues of $g$
and let $(u,u')$ be the corresponding eigenvectors,
so that $gu=ru$, $gu'=r^{-1}u'$. Let $u^*:\R^2\ra\R$
be the functional defined by $u^*(u)=1$, $u^*(u')=0$. 
We can choose $u$ in such a way that $u^*(v_0)>0$, where
$v_0=v_{\FF}$ (note that $v_0$ cannot be an eigenvector of
$g$ since $M=\chi(v_0,g v_0)>0$).
Consider the half-plane $H=\{v\in\R^2|\ u^*(v)>0\}$.
Let $\CC=\CC^{\th}[m]$
for appropriate $\th\in\R\cup\{\infty\}$ and $m\in\Z$, so that 
$\FF\in\CC$ and the vectors $(v_{\GG}, \GG\in\CC)$
belong to the closure of $H$. Then $F$ preserves $\CC$ (see 
Proposition \ref{eigenprop}).
For an object $\GG$ of $D^b(X)$ we set 
$$\rk_{\CC}(\GG)=\frac{u^*(v_{\GG})}{u^*(v_0)},$$
so that $\rk_{\CC}(\GG)\ge 0$ for all $\GG\in\CC$.
From the definition of $u^*$ we immediately derive that
$$\rk_{\CC}(F(\GG))=r\cdot\rk_{\CC}(\GG).$$
Note that $\CC$ contains the subcategory equivalent
to the category of stable bundles on a noncommutative $2$-torus and 
$\rk_{\CC}$ is proportional to the rank function on such bundles.
Let us denote $\FF_n=F^n(\FF)\in\CC$. Since $\rk_{\CC}(\FF)=1$, we obtain
that
$$\rk_{\CC}(\FF_n)=r^n.$$

To prove the above theorem we are going to use
the twist functors $R_{\FF_n}:D^b(X)\ra D^b(X)$. More precisely,
let us consider the objects 
$$\FF'_n:=R_{\FF_n}R_{\FF_{n-1}}\ldots R_{\FF_1}(\FF_0)\in D^b(X),$$
where $n>0$.
It is convenient to extend this definition to $n=0$ by
setting $\FF'_0=\FF_0$. As we will see, the properties of the algebra
$A$ depend on whether some (or all) $\FF'_n$
belong to the subcategory $\CC$ and also on the vanishing of
some (or all) spaces $\Hom^1(\FF'_n,\FF_m)$ for $n<m$.

\begin{lem}\label{Hilblem}
Consider the following generating series
$$F(t,u)=\sum_{n\ge 0,k\ge 0}\chi(\FF'_n,\FF_{n+k+1})t^nu^k,$$
$$R(t)=\sum_{n\ge 0}\rk_{\CC}(\FF'_n)t^n.$$
Then one has
$$F(t,u)=\frac{M(1+tu)}{(1-Nu+u^2)(1-(M-N)t+t^2)},$$
$$R(t)=\frac{1+r^2t}{1-(M-N)rt+r^2t^2}.$$
\end{lem}

\Pf .
By the definition of $\FF'_n$ we
have an exact triangle
\begin{equation}\label{extriangle}
\FF'_{n-1}\ra
\Hom^{\bullet}(\FF'_{n-1},\FF_n)^*\otimes \FF_n\ra \FF'_n\ra \FF'_{n-1}[1].
\end{equation}
This implies the following relations:
$$\chi(\FF'_n,\FF_m)=\chi(\FF'_{n-1},\FF_n)\chi(\FF_n,\FF_m)-
\chi(\FF'_{n-1},\FF_m),\ m>n\ge 1,$$
$$\rk_{\CC}(\FF'_n)=\chi(\FF'_{n-1},\FF_n)r^n-\rk_{\CC}(\FF'_{n-1}), n\ge 1.$$
Note that 
$\chi(\FF_n,\FF_m)=\chi(\FF_0,\FF_{m-n})$ for $m>n$
is a coefficient of the Hilbert series of $A$:
$$H(t)=1+\sum_{n\ge 1}\chi(\FF_0,\FF_n)t^n.$$
Therefore, denoting
$$F_0(t)=F(t,0)=\sum_{n\ge 0}\chi(\FF'_n,\FF_{n+1})t^n,$$
we obtain the equations
$$(u+t)F(t,u)=H(u)(1+tF_0(t))-1,$$
$$R(t)=1+rtF_0(rt)-tR(t).$$
Substituting $u=-t$ into the first equation we get
$$F_0(t)=\frac{H(-t)^{-1}-1}{t},$$
and therefore,
$$F(t,u)=\frac{H(u)H(-t)^{-1}-1}{u+t},$$
$$R(t)=\frac{H(-rt)^{-1}}{1+t}.$$
It remains to use the formula
$$H(t)=\frac{1+(M-N)t+t^2}{1-Nt+t^2}$$
that was proven in Proposition \ref{eigenprop}.
\ed

\noindent
{\it Proof of Theorem \ref{Koszulthm}.}
(a) If the map $A_1\otimes A_1\ra A_2$ is
surjective then $M^2=(\dim A_1)^2\ge\dim A_2=MN$, hence
$M\ge N$.

Conversely, assume that $M\ge N$.
Then we claim that $\FF'_1\in\CC$ and 
$\Hom^1(\FF'_1,\FF_m)=0$ for $m>2$.
Indeed, since $\FF'_1=R_{\FF_1}(\FF_0)$, it is a stable object
of $D^b(X)$. Hence, the exact triangle (\ref{extriangle})
for $n=1$ implies that either
$\FF'_1\in\CC$ or $\FF'_1\in\CC[1]$. To check that $\FF'_1\in\CC$
it suffices to prove the inequality $\rk_{\CC}(\FF'_1)>0$.
But
$$\rk_{\CC}(\FF'_1)=Mr-1=(M-N)r+r^2\ge r^2>0$$
which proves our first claim
(we used the equality $N=r+r^{-1}$).
On the other hand, since $\FF'_1$ and $\FF_m$ are both stable objects of $\CC$,
the vanishing of $\Hom^1(\FF'_1,\FF_m)$ would follow
from the inequality $\chi(\FF'_1,\FF_m)>0$ (see Lemma \ref{nzhomlem}).
From Lemma \ref{Hilblem} we get
$$\sum_{k\ge 0}\chi(\FF'_1,\FF_{k+2})u^k=
\frac{\partial F}{\partial t}(0,u)=\frac{M(M-N+u)}{1-Nu+u^2}.$$
The latter series has positive coefficients except maybe for
the constant term. Hence, $\Hom^1(\FF'_1,\FF_m)=0$ for $m>2$.
Now from the exact sequence
\begin{equation}\label{exactseq}
0\ra\Hom^0(\FF'_1,\FF_m)\ra\Hom^0(\FF_0,\FF_1)\otimes\Hom^0(\FF_1,\FF_m)\ra
\Hom^0(\FF_0,\FF_m)\ra\Hom^1(\FF'_1,\FF_m)
\end{equation}
we derive the surjectivity of the map $A_1\otimes A_{m-1}\ra A_m$
for $m>2$. In the case $M>N$ the above argument
shows that $\Hom^1(\FF'_1,\FF_2)=0$, hence in this case the map
$A_1\otimes A_1\ra A_2$ is also surjective.
On the other hand, if $M=N$ then $v_{\FF'_1}=v_{\FF_2}$, so
either $\Hom^{\bullet}(\FF'_1,\FF_2)=0$, or
$\FF'_1\simeq \FF_2$. Hence, in this case the map $A_1\otimes A_1\ra A_2$
is surjective if and only if $\det(\FF'_1)\not\simeq \det(\FF_2)$.
Using the triangle (\ref{extriangle}) we get that
$\det(\FF'_1)\simeq\det(\FF_1)^N\otimes\det(\FF_0)^{-1}$ which
leads to the condition in the formulation of part (a).

\noindent
(b) By the result of part (a) we can assume that $A$ is generated
by $A_1$. The statement
that the algebra $A$ is quadratic is equivalent to surjectivity
of the natural maps
\begin{equation}\label{quadrmaps}
A_{m-2}\otimes I\ra\ker (A_{m-1}\otimes A_1\ra A_m)
\end{equation}
for all $m\ge 3$, where $I=\ker(A_1\otimes A_1\ra A_2)$ is
the space of quadratic relations.
From the exact sequences (\ref{exactseq}) above we see
that $I$ can be identified with $\Hom^0(\FF'_1,\FF_2)$ and
that the map (\ref{quadrmaps}) can be identified with the natural map
$$\Hom^0(\FF_2,\FF_m)\otimes\Hom^0(\FF'_1,\FF_2)\ra
\Hom^0(\FF'_1,\FF_m),\ m\ge 3.$$
Now the exact triangle (\ref{extriangle}) for $n=2$
shows that the kernel and the cokernel of this map are
$\Hom^0(\FF'_2,\FF_m)$ and $\Hom^1(\FF'_2,\FF_m)$ respectively
(we use the facts about $\FF'_1$ proven in part (a)).
Therefore, to prove that the algebra $A$ is quadratic it suffices to
show that $\Hom^1(\FF'_2,\FF_m)=0$ for $m\ge 3$.

The same argument as in (a) shows that either $\FF'_2\in\CC$ or
$\FF'_2\in\CC[1]$.  Moreover, we have
$$\rk_{\CC}(\FF'_2)=((M-N)^2-1)r^2+(M-N)r^3,$$
so for $M\ge N+1$ we have $\rk_{\CC}(\FF'_2)>0$, hence
$\FF'_2\in\CC$. On the other hand, using Lemma \ref{Hilblem} we get
$$\sum_{k\ge 0}\chi(\FF'_2,\FF_{k+3})u^k=
\frac{(M-N)^2-1+(M-N)u}{1-Nu+u^2}.$$
Thus, if $M-N\ge 2$ then $\chi(\FF'_2,\FF_m)>0$ for all $m\ge 3$.
Since $\FF'_2$ is a stable object of $\CC$, this implies
the required vanishing of $\Hom^1(\FF'_2,\FF_m)$. If $M=N+1$ then
$\chi(\FF'_2,\FF_m)>0$ for $m\ge 4$ while $\chi(\FF'_2,\FF_3)=0$.
Hence, in this case $\Hom^1(\FF'_2,\FF_m)=0$ for $m\ge 4$
and either $\Hom^{\bullet}(\FF'_2,\FF_3)=0$ or $\FF'_2\simeq \FF_3$.
The latter isomorphism occurs exactly when $\det(\FF'_2)\simeq\det(\FF_3)$.
It remains to use (\ref{extriangle}) to get an isomorphism
$$\det(\FF'_2)\simeq\det(\FF_2)^{N+1}\otimes\det(\FF'_1)^{-1}\simeq
\det(\FF_2)^{N+1}\otimes\det(\FF_1)^{-N-1}\otimes\det(\FF_0).$$

Conversely, if the algebra $A$ is quadratic then
$$\chi(\FF'_2,\FF_3)=M\chi(\FF'_1,\FF_2)-\chi(\FF'_1,\FF_3)\ge 0,$$
hence $(M-N)^2\ge 1$. Since we already know that $M-N\ge 0$
this implies that $M-N\ge 1$.

\noindent
(c) If the algebra $A$ is Koszul then its Hilbert series $H(t)$
has the property that $H(-t)^{-1}$ has non-negative coefficients.
But
$$H(-t)^{-1}=\frac{1+Nt+t^2}{1-(M-N)t+t^2}=1+\frac{Mt}{1-(M-N)t+t^2},$$
so this series has non-negative coefficients only if $M-N\ge 2$.

Conversely, assume that $M-N\ge 2$. We claim that all the objects $\FF'_n$
belong to $\CC$ and that $\Hom^1(\FF'_n,\FF_m)=0$ 
for $m>n$. We argue by induction.
Assume that the assertion is true for $n-1$. Looking
at the exact triangle (\ref{extriangle}) we conclude as in part (a) that
either $\FF'_n\in\CC$ or $\FF'_n\in\CC[1]$.
Since by Lemma \ref{Hilblem} we also have $\rk_{\CC}(\FF'_n)>0$
this implies that $\FF'_n\in\CC$.
On the other hand, the same Lemma shows that
$\chi(\FF'_n,\FF_m)>0$ for $m>n$. By Lemma \ref{nzhomlem}
this implies the vanishing of $\Hom^1(\FF'_n,\FF_m)$ for $m>n$.

Now let us set
$$K_n:=\oplus_{m>n}\Hom^0(\FF'_n,\FF_m).$$
Then $K_n$ has a natural structure of (graded) left $A$-module and
from (\ref{extriangle}) we deduce the following exact sequences
of $A$-modules:
$$0\ra K_n\ra A(-n)\otimes\Hom(\FF'_{n-1},\FF_n)\ra K_{n-1}\ra 0,\ n\ge 1,$$
where $A(-n)$ is the free $A$-module with one generator in degree $n$.
Since $K_0$ is the augmentation ideal $A_{+}\sub A$, putting
these sequences together we obtain a free resolution of the trivial
module $k$ of the form
$$\ldots\ra A(-n)\otimes\Hom(\FF'_{n-1},\FF_n)\ra\ldots\ra
A(-1)\otimes\Hom(\FF'_0,\FF_1)\ra A.$$
This implies that $A$ is Koszul and 
$A^!_n\simeq\Hom(\FF'_{n-1},\FF_n)^*$ for $n\ge 1$.

To prove the last statement of the theorem we observe that the
exact triangle (\ref{extriangle}) shows that
$$A^!_n\simeq\Hom(\FF'_{n-1},\FF_n)^*\simeq\Hom(\FF_n,\FF'_n).$$
Since $R_{\FF_n}\simeq F^n R_{\FF} F^{-n}$, we also have
$$\FF'_n\simeq F^n (R_{\FF} F^{-1})^n(\FF).$$
Therefore,
$$A^!_n\simeq\Hom(\FF_n,\FF'_n)\simeq \Hom(\FF,(R_{\FF}F^{-1})^n(\FF)).$$
It is easy to check that this isomorphism is compatible with
the multiplication on $A^!$ and on $A_{R_{\FF}F^{-1},\FF}$.
\ed

\begin{exs} 1. If $\a_F=1$ then by Proposition \ref{cubefunprop}
part (b) of the above theorem states that $A_{F,\FF}$ is quadratic iff 
$M\ge N+2$. For example,
this is the case for $F=(\otimes L)$, where $L$ is a line bundle.
This leads to the well-known statement
that the corresponding algebra
$A_{F,\OO_X}=\oplus_{n\ge 0}H^0(X,L^n)$ 
is quadratic iff $\deg(L)\ge 4$ (then it is also
Koszul). More generally, if $F=(\otimes L)\circ\si^*$, where $\si$
is an automorphism of $X$, then $A_{F,\OO_X}$ is the so called twisted
coordinate algebra attached to the pair $(L,\si)$. Such algebras were
considered in \cite{AVdB}. For example, if $\si$ is a translation
then such an algebra is quadratic iff $\deg(L)\ge 4$,
in which case it is also Koszul.

2. Consider the case $M=N+1$, $\a_F=-1$. Then the algebra $A$ is
often quadratic but never Koszul. The quadratic dual has $A^!_3=0$.
For example, if $L$ is a line bundle of degree $3$ such that
$[-1]^*L\otimes L^{-1}$ is not of order $2$, then these conditions
are satisfied for $F=(\otimes L)\circ [-1]^*$. The corresponding algebra is
$$A=k\oplus H^0(L)\oplus H^0(L\otimes [-1]^*L)\oplus
H^0(L\otimes [-1]^*L\otimes L)\oplus\ldots$$
with the multiplication rule $f*g=f\cdot [(-1)^{\deg(f)}]^*g$.

3. If $M=N+1$ and $\a_F=1$ then the algebra $A$ is not quadratic:
one has to add one cubic relation to the quadratic relations.
\end{exs}

The following result allows to check under what conditions
Theorem \ref{Koszulthm} can be applied to sufficiently
high powers of a given autoequivalence.

\begin{prop}\label{Verprop} Assume that an element 
$g\in\SL_2(\Z)$ has positive real eigenvalues and
that the vector $v\in\Z^2\setminus\{0\}$ satisfies
$M:=\chi(v,gv)>0$. 
Then the following conditions are equivalent:

\noindent
(i) $\chi(v,g^n v)-\tr(g^n)\to +\infty$ as $n\to +\infty$;

\noindent
(ii) $\chi(v,g^n v)-\tr(g^n)\ge 0$ for some $n>0$;

\noindent
(iii) $M>r_1-r_2$, where $r_1\ge r_2$ are eigenvalues of $g$;

\noindent
(iv) either $g$ is unipotent, or $M\ge N$, where $N:=\tr(g)$.
\end{prop}

\Pf .  Note that by Lemma \ref{eigenlem} we have
$\chi(v,g^n v)>0$ for all $n>0$. Moreover,
from the proof of that Lemma we get
$$\sum_{n\ge 1} \chi(v,g^n v) t^n=\frac{Mt}{1-Nt+t^2}.$$

If $g$ is unipotent then we have $N=\tr(g^n)=2$, while
$\chi(v,g^n v)=nM$, so the assertion is clear. 

Now assume that $g$ has two distinct
eigenvalues $r>r^{-1}>0$. Then from the above formula we get
$$\chi(v,g^n v)=M\cdot\frac{r^n-r^{-n}}{r-r^{-1}}.$$
On the other hand, $\tr(g^n)=r^n+r^{-n}$.
Using these formulas it is not difficult to show
the equivalence of (i)--(iv).
Indeed, clearly, (i) implies (ii). The implication (ii)$\implies$(iii)
follows immediately from the chain of inequalities
$$\frac{Mr^n}{r-r^{-1}}>\frac{M(r^n-r^{-n})}{r-r^{-1}}\ge r^n+r^{-n}>r^n.$$
To prove (iii)$\implies$(iv) we note that
$r-r^{-1}=\sqrt{N^2-4}$. Thus, the inequality $M>r-r^{-1}$ implies that
$M^2>N^2-4$, hence $M\ge N$ (since $N\ge 3$ in our case).
Finally, if $M\ge N$ then
$M>r-r^{-1}$, in which case
$\chi(v,g^n v)-\tr(g^n)\to +\infty$ as $n\to +\infty$.
This proves (iv)$\implies$(i).
\ed

\section{Ampleness and noncommutative $\Proj$}\label{projsec}

\subsection{Ampleness criterion}

Let $A$ be a graded $k$-algebra of the form $A=\oplus_{i\ge 0}A_i$, 
where $A_0=k$. Recall that a finitely generated right graded $A$-module
$M$ is called {\it coherent}
if for every finite collection of homogeneous elements 
$m_1,\ldots,m_n\in M$ the (right) $A$-module of relations between 
$m_1,\ldots,m_n$ is finitely generated. Coherent modules form
an abelian subcategory in the category of all modules.
An algebra $A$ is called
right coherent if it is finitely generated and is coherent as a 
right module over itself.
We denote by $\cohproj A$ the quotient of the category of coherent
$A$-modules by the Serre subcategory of bounded coherent modules.
Below we are going to show that the categories $\CC^{\th}$,
where $\th$ is a quadratic irrationality (or a rational number),
are equivalent to such quotient categories for appropriate 
algebras of the form $A_{F,\FF}$.

Let $\CC$ be an abelian category equipped with an autoequivalence
$F:\CC\ra\CC$. For simplicity we will assume that $F$ is an automorphism
of $\CC$
(the general case can be reduced to this one, see \cite{AZ}).
Under appropriate ampleness assumptions
the category $\CC$ can be recovered
from the algebra $A_{F,O}$, where $O$ is an object of $\CC$.
Recall (see \cite{SVdB}, \cite{P-coh}) that a sequence
$(O_n,n\in\Z)$ of objects of $\CC$ is called {\it ample} if the following
two conditions hold: (i) for every surjection $X\ra Y$ in $\CC$
the induced map $\Hom_{\CC}(O_n,X)\ra\Hom_{\CC}(O_n,Y)$ is surjective
for all $n<<0$; (ii) for every object $X\in\CC$ and every $n\in\Z$
there exists a surjection $\oplus_{j=1}^s O_{i_j}\ra X$, where $i_j<n$
for all $j$. The main theorem of \cite{P-coh} implies that if the 
sequence $(F^n O,n\in\Z)$ is ample then the algebra $A_{F,O}$
is right coherent and the categories $\CC$ and $\cohproj A_{F,O}$ are equivalent.
Similar result for Noetherian categories was proven by Artin and Zhang
in \cite{AZ}. The following proposition shows that 
we do need a more general theorem of \cite{P-coh}, since 
for irrational $\th$ the
categories $\CC^{\th}$ are not Noetherian.

\begin{prop}\label{noethprop} Assume that $\th$ is irrational.
Then every non-zero object in $\CC^{\th}$
is not Noetherian.
\end{prop}

\Pf . It suffices to prove that every stable object $\FF\in\CC^{\th}$ 
is not Noetherian. Recall that the vector 
$v_{\FF}=(\deg(\FF),\rk(\FF))\in\Z^2$ satisfies
$\deg(\FF)-\rk(\FF)\th>0$. Moreover, since $\FF$ is stable, the numbers 
$\deg(\FF)$ and $\rk \FF$ are relatively prime. Let $(m,n)$ be the unique pair
of integers such that $m\rk(\FF)-n\deg(\FF)=1$ and
$0<m-n\th<\deg(\FF)-\rk(\FF)$. There exists a stable object $\FF'\in\CC^{\th}$
with $v_{\FF'}=(m,n)$ (see Proposition \ref{cohdim}). 
By Lemma \ref{nzhomlem} one has $\dim\Hom(\FF,\FF')=1$ and
$\Hom^1(\FF,\FF')=0$. This implies that there is an exact triangle
$$\FF'[-1]\ra L_{\FF}(\FF')\ra \FF\ra \FF'\ra\ldots$$
where $L_{\FF}$ is the left twist functor corresponding to $\FF$. 
Note that the object $L_{\FF}(\FF')$ is stable, so
either $L_{\FF}(\FF')\in\CC^{\th}$, or $L_{\FF}(\FF')\in\CC^{\th}[-1]$.
But the vector $v_{L_{\FF}(\FF')}=v_{\FF}-v_{\FF}'$ lies in the half-plane
$\{(x,y)|\ x-y\th>0\}$, hence $L_{\FF}(\FF')$ is in $\CC^{\th}$.
This means that $\FF'$ is a proper quotient-object of $\FF$ in $\CC^{\th}$.
Iterating this procedure we will obtain an infinite sequence
$\FF\ra \FF_1\ra \FF_2\ra\ldots$, where $\FF_{n+1}$ is a proper quotient
of $\FF_n$.
\ed

The following theorem gives a criterion of ampleness for sequences
of the form $(F^n \FF, n\in\Z)$ in the categories $\CC^{\th}$.

\begin{thm}\label{amplethm} 
Let $F:D^b(X)\ra D^b(X)$ be an autoequivalence
such that the element $g=\pi(F)\in\SL_2(\Z)$
has distinct positive real eigenvalues.
Let $u=(x,y)\in\R^2$ be an eigenvector of $g$
with the eigenvalue $<1$ and let $\th=x/y$
(if $y=0$ then $\th=\infty$). Let $\FF_0$ be a stable object
of $\CC^{\th}$ and let $v_0=(\deg \FF_0,\rk \FF_0)$ be the corresponding
primitive vector in $\Z^2$. 
Assume that $F(\FF_0)\in\CC^{\th}$.
Denote also $N=\tr(g)$ and $M=\chi(\FF_0,F(\FF_0))=\chi(v_0,gv_0)$.

\noindent
(a) If $M\ge N-1$ then the sequence
$(F^n(\FF_0),n\in\Z)$ in $\CC^{\th}$ is ample.

\noindent
(b) If $0<M<N-1$ then the algebra $A_{F,\FF_0}$ is not finitely
generated, hence, the sequence $(F^n(\FF_0))$ is not ample. 
\end{thm}

\Pf . Below we will denote the coordinates of a vector 
$v\in\R^2$ by $(\deg(v),\rk(v))$. 
Let us denote $\FF_n=F^n(\FF_0)$, $v_n=v_{\FF_n}=g^n v_0$.
By Proposition \ref{tstrprop} one has $F(\CC^{\th})=\CC^{\th}[m]$
for some $m\in\Z$. Hence, our assumption $F(\FF_0)\in\CC^{\th}$ implies
that $F(\CC^{\th})=\CC^{\th}$. In particular, $\FF_n\in\CC^{\th}$
for all $n\in\Z$. 

\noindent
(a) Assume first that $M\ge N$. Note that
$v_0$ is not an eigenvector of $g$, so we can choose $u$
in such a way that $\chi(u,v_0)>0$. Then $\chi(u,v_n)>0$ for all $n\in\Z$.
Moreover, since $\chi(g^{-1}v_0,v_0)>0$ and since $u$ is an eigenvector
of $g^{-1}$ with the eigenvalue $>1$, it follows that
$\deg(v_n)/\rk(v_n)$ tends to $\th$ as $n\to-\infty$.
Observe also that all the vectors
$\{v_{\FF}, \FF\in\CC^{\th}\}$ lie in the half-plane 
$H=\{v:\ \chi(u,v)>0\}$ (since $\th$ is irrational). 
It suffices to prove that for every $\FF\in\CC^{\th}$ the following holds:

\noindent
(i) $\Hom^1(\FF_n,\FF)=0$ for $n<<0$;

\noindent
(ii) the natural map $\Hom(\FF_{n},\FF)\otimes \FF_{n}\ra \FF$ is surjective
for $n<<0$.

Moreover, it is enough to check these statements for stable $\FF$. Also,
for (ii) it is enough to prove that $L_{\FF_{n}}(\FF)$ is in $\CC^{\th}$ for
$n<<0$,
where $L_{\FF_n}$ is the left twist functor associated with $\FF_n$.
Since the vectors $v_n$ lie in the half-plane $H$,
$\chi(v_{n-1},v_n)>0$ and $\deg(v_n)/\rk(v_n)\to\th$ as $n\to-\infty$,
it follows that for every vector $v$ in $H$
one has $\chi(v_n,v)>0$ for $n<<0$. Applying Lemma \ref{nzhomlem}
to $\FF_n$ and a stable object $\FF\in\CC^{\th}$ we immediately derive (i).
It remains to prove that for such $\FF$ one has
$L_{\FF_{n}}(\FF)\in\CC^{\th}$ for $n<<0$. 
Since $L_{\FF_n}(\FF)$ is a stable object that fits into an exact
triangle
$$\FF[-1]\ra L_{\FF_n}(\FF)\ra \Hom(\FF_n,\FF)\otimes \FF_n\ra \FF,$$
it suffices to prove that $v_{L_{\FF_n}(\FF)}$ belongs to $H$. But
$$v_{L_{\FF_n}(\FF)}=\chi(\FF_n,\FF) v_n-v_{\FF}=
\chi(v_n,v_{\FF})v_n-v_{\FF},$$
so
$$\chi(u,v_{L_{\FF_n}(\FF)})=\chi(v_n,v)\chi(u,v_n)-\chi(u,v).$$
where $v=v_{\FF}$.

Let $r<1$ be the eigenvalue of $g$
corresponding to $u$ and let $u'$ be the eigenvector corresponding
to $r^{-1}$. Rescaling $u$ and $u'$ we can assume
that $v_0=u+u'$. Then $v_n=r^nu+r^{-n}u'$,
$\chi(u,v_n)=r^{-n}\chi(u,u')>0$ and
$$\chi(v_n,v)=r^n\chi(u,v)+r^{-n}\chi(u',v).$$
It follows that 
$$\chi(u,v_{L_{\FF_n}(\FF)})=
(\De-1)\chi(u,v)+r^{-2n}\De\chi(u',v),$$
where $\De=\chi(u,u')$.
Since $\chi(u,v)>0$, this quantity is positive for $n<<0$ provided that
$\De>1$. But
$$\De=\frac{\chi(v_0,gv_0)}{r^{-1}-r}\ge\frac{r+r^{-1}}{r^{-1}-r}>1,$$
where we used our assumption $M\ge N$.

Now let us consider the case $M=N-1$. The condition (i)
is still satisfied, however (ii) has to be replaced by a weaker condition.
For every $\FF\in\CC^{\th}$ and $n\in\Z$ let us set
$$T_n\FF:=\coker(\Hom(\FF_n,\FF)\otimes \FF_n\ra\FF).$$
It suffices to prove that for every $\FF$ we have
$$T_{-n+m}\ldots T_{-1+m}T_{m}\FF=0$$
for some $m\in\Z$ and some $n>0$.
As before we can assume that
$\FF$ is a stable object in $\CC^{\th}$.
Using the action of $F$ we can reduce ourselves to the case
when the vector $v=v_{\FF}$ satisfies $\chi(v_0,v)>0$.
In this case we will show that
$T_{-n}\ldots T_{-1}T_0\FF=0$ for some $n>0$.
By Lemma \ref{nzhomlem}, we have $\Hom^1(\FF_0,\FF)=0$,
so there is an exact triangle
$$L_{\FF_0}(\FF)\ra \Hom(\FF_0,\FF)\otimes \FF_0\ra\FF\ra 
L_{\FF_0}(\FF)[1].$$
Let $h_i\in\SL_2(\Z)$ be the matrix corresponding to the functor 
$L_{\FF_i}[1]:D^b(X)\ra D^b(X)$.
We claim that if $v_{L_{\FF_0}(\FF)[1]}=h_0(v)\not\in H$ then $T_0\FF=0$.
Indeed, since $L_{\FF_0}(\FF)$ is a stable object, this would imply
that $L_{\FF_0}(\FF)\in\CC^{\th}$, so the map $\Hom(\FF_0,\FF)\otimes \FF_0\ra\FF$
is surjective. Otherwise, we have $v_{L_{\FF_0}(\FF)[1]}\in H$
which implies that $L_{\FF_0}(\FF)[1]$ belongs to $\CC^{\th}$, the 
map $\Hom(\FF_0,\FF)\otimes\FF_0\ra\FF$ is injective with the cokernel
$$T_0\FF\simeq L_{\FF_0}(\FF)[1]$$
and $v_{T_0\FF}=h_0(v)$.
Continuing to argue in this way we see that it is enough to show
the existence of $n>0$ such that
$h_{-n}\ldots h_{-1}h_0(v)\not\in H$.
Using the formula
$v_{L_{\FF_0}(\FF)[1]}=v-\chi(v_0,v)v_0$
we can write the matrix of $h_0$ with respect to 
the basis $(u,u')$:
$$h_0=\left(\matrix 1+\De & -\De \\ \De & 1-\De\endmatrix\right),$$
where $\De=\chi(u,u')$. 
Similarly,
$$h_{-i}=\left(\matrix r^{-i} & 0 \\ 0 & r^i\endmatrix\right)h_0
\left(\matrix r^i & 0 \\ 0 & r^{-i}\endmatrix\right).$$
Therefore,
$$h_{-n}\ldots h_{-1}h_0=
\left(\matrix r^{-n+1} & 0 \\ 0 & r^{n-1}\endmatrix\right)S^{n+1},$$
where
$$S=\left(\matrix r & 0 \\ 0 & r^{-1}\endmatrix\right) h_0=
\left(\matrix r(1+\De) & -r\De \\ r^{-1}\De & r^{-1}(1-\De)\endmatrix\right).$$
But $\det(S)=1$ and
$$\tr(S)=r(1+\De)+r^{-1}(1-\De)=N+(r-r^{-1})\De=N-M=1.$$
Hence, $S^2-S+1=0$ and therefore $S^3=-1$.
It follows that $h_{-2}h_{-1}h_0(v)$ is not in $H$ which finishes the proof.

\noindent
(b) Let $A=A_{F,\FF_0}$. For a graded right $A$-module 
$M$ and $n\in\Z$ let us set
$$T_nM:=\coker(M_n\otimes A(-n)\ra M),$$
where $A(-n)$ is a free $A$-module with $A(-n)_i=A_{i-n}$.
To show that the algebra $A$ is not finitely generated
we have to show that for all $n\ge 1$ one has
$$T_nT_{n-1}\ldots T_1 A_{\ge 1}\neq 0,$$
where $A_{\ge 1}=\oplus_{i\ge 1}A_i$.

For every $\FF\in\CC^{\th}$ and $n\in\Z$ we set
$$\Ga_{\ge n}(\FF)=\oplus_{m\ge n}\Hom(\FF_{-m},\FF).$$
This space has a natural structure of a graded right $A$-module.
For example, we have $\Ga_{\ge n}(\FF_{-n})=A(-n)$.
Now we claim that it is enough to prove that 
$h_{-n}h_{-n+1}\ldots h_{-1}(v_0)$ is a nonzero vector in
$H$ for all $n\ge 1$,
where we use the notation from the proof of part (a).
Indeed, as we have seen above this would imply that for every $n\ge 1$
the object $T_{-n}T_{-n+1}\ldots T_1\FF_0$ is stable and that
we have exact sequences
$$0\ra\Hom(\FF_{-n},T_{-n+1}\ldots T_1\FF_0)\otimes \FF_{-n}\ra
T_{-n+1}\ldots T_1\FF_0\ra T_{-n}T_{-n+1}\ldots T_1\FF_0\ra 0$$
in $\CC^{\th}$. Using these exact sequences and Lemma \ref{nzhomlem}
one can easily see that 
$$\Ga_{\ge n+1}(T_{-n}T_{-n+1}\ldots T_1\FF_0)\simeq
T_nT_{n-1}\ldots T_1 A_{\ge 1}.$$
Again applying Lemma \ref{nzhomlem} we conclude that this space is
not zero for every $n\ge 1$ which proves our claim.

To prove that $h_{-n}\ldots h_{-1}(v_0)$ is a nonzero vector in $H$
it suffices to show that $S^n(\R_{>0}u+\R_{>0}u')\sub H$
for all $n\ge 1$. Note that $\tr(S)=N-M\ge 2$, so $S$
has real positive eigenvalues.
Since $\chi(u,Su)=r^{-1}\De^2>0$, it is enough to prove that there exists
an eigenvector of $S$ of the form $-xu+u'$, where $x>0$.
Equivalently, the equation
$$\chi(-xu+u',S(-xu+u'))=r^{-1}\De x^2+[r(1+\De)-r^{-1}(1-\De)]x+r\De=0$$
should have a positive root. For this two inequalities should hold:
$D=b^2-4\De^2\ge 0$ and $b<0$, where $b=r(1+\De)-r^{-1}(1-\De)$.
But
$$b=r-r^{-1}+N\De=\frac{NM}{r^{-1}-r}-r^{-1}+r=\frac{NM-N^2+4}{r^{-1}-r}<0$$
since $M\le N-2$ and $N>2$. Finally,
$$D=(r-r^{-1}+N\De)^2-4\De^2=N^2-4-2NM+(N^2-4)\De^2=N^2-4-2NM+M^2\ge 0$$
since $N-M\ge 2$.
\ed

\begin{rem} If the equivalent conditions of the above theorem 
are satisfied then we also have $F^n(\FF_0)\in\CC^{\th'}$
for all $n$, where $(\th',1)$ is the eigenvector of $g$ 
corresponding to the eigenvalue
$>1$. The sequence $(F^n(\FF_0))$ is not ample in $\CC^{\th'}$, 
since there are objects 
$\FF\in\CC^{\th'}$ with $\Hom^1(F^n(\FF_0),\FF)\neq 0$ and 
$\Hom(F^n(\FF_0),\FF)=0$ for all $n<<0$. Nevertheless,
we still have an equivalence of the derived category of $\CC^{\th'}$
with the derived category of $\cohproj A_{F,\FF}$, since
both categories are equivalent to $D^b(X)$. It would be interesting
to find a general framework for this kind of equivalences associated
with non-ample sequences.
\end{rem}

\begin{cor} Let $(F,\FF)$ be a pair 
satisfying the equivalent conditions of Proposition \ref{eigenprop}
and let $\pi(F)=g\in\SL_2(\Z)$, $N=\tr(g)$,
$M=\chi(\FF,F(\FF))$. Then the algebra $A_{F,\FF}$ is finitely
generated if and only if $M\ge N-1$.
\end{cor}

\Pf . If $g$ has distinct eigenvalues then this
follows from Theorem \ref{amplethm}. 
Now assume that $g$ is unipotent (so that $N=2$).
Then the statement reduces to the case when $F$ is a composition
of the tensoring by a line bundle $L$ with an automorphism
of $X$. In this case we can assume that $\FF\in\Coh X$.
Since $M\ge 1$ it follows that $\FF$ is a vector bundle and
$\deg(L)\ge 1$. It easy to see that in this case
the sequence $(F^n(\FF))$ is ample, hence, the algebra
$A_{F,\FF}$ is finitely generated.
\ed

\subsection{Projectivity of $\CC^{\th}$}

Now we can show that every category $\CC^{\th}$,
where $\th$ is a quadratic irrationality, can be described
as a ``noncommutative $\Proj$''.

\begin{thm}\label{projthm1} 
For every quadratic irrationality $\th\in\R$
there exists 
an autoequivalence $F:D^b(X)\ra D^b(X)$ preserving $\CC^{\th}$ and a
stable
object $\FF\in\CC^{\th}$ such that the sequence $(F^n\FF, n\in\Z)$ is ample.
Hence, the corresponding algebra $A_{F,\FF}$ is right coherent
and $\CC^{\th}\simeq\cohproj A_{F,\FF}$.
\end{thm}

\Pf . Let $\a\th^2+\b\th+\ga=0$ be the equation satisfied by
$\th$, where $\a,\b,\ga\in\Z$, $\a>0$. Consider the ring 
$R=\Z[\a\th]\sub\Q(\th)$. Then $R$ is contained in $\Z+\Z\th$
and $R(\Z+\Z\th)\sub\Z+\Z\th$. Let $r\in R^*$ be a unit such that
$0<r<1$ (such a unit always exists). Then the multiplication by $r$ induces
an invertible operator on $\Z+\Z\th$ with determinant equal to
$\Nm(r)=1$. Hence, we have $r=c\th+d$, $r\th=a\th+b$ for some 
$g=\left(\matrix a & b \\ c & d \endmatrix\right)\in\SL_2(\Z)$.
Then $u=(\th,1)\in\R^2$ is an eigenvector of $g$ corresponding to the
eigenvalue $r$. We claim that there exists
a primitive vector $v\in\Z^2$ such that $\chi(u,v)>0$ and 
$\chi(v,gv)\ge\tr(g)$. Indeed, let 
$u'\in\R^2$ be an eigenvector of $g$ corresponding to the eigenvalue $r^{-1}$
and such that $\chi(u,u')>0$. We can find a primitive
vector $v\in\Z^2$ such that $v=xu+yu'$, where $x>0$, $y>0$, and 
$xy\ge\frac{r+r^{-1}}{(r^{-1}-r)\chi(u,u')}$. Then
$\chi(u,v)=y\chi(u,u')>0$ and
$$\chi(v,gv)=\chi(xu+yu',xru+yr^{-1}u')=xy(r^{-1}-r)\chi(u,u')\ge r+r^{-1}$$
as required.
It remains to choose an autoequivalence $F$ with $\pi(F)=g$ such that
$F$ preserves $\CC^{\th}$ and an object $\FF\in\CC^{\th}$ 
with $v_{\FF}=v$, and then
apply Theorem \ref{amplethm}.
\ed 

Finally, we are going to
show that the category $\CC^{\th}$ 
for arbitrary $\th\in\R$ can be
represented in the form $\cohproj A$ for some coherent
$\Z$-algebra $A$. Recall that the notion of $\Z$-algebra 
is a natural generalization of the notion of graded algebra
(see \cite{BP}, \cite{P-coh}): such an algebra is equipped
with a decomposition $A=\oplus_{i\le j}A_{i,j}$ and the case
of a graded algebra corresponds to $A_{i,j}=A_{j-i}$.
As in the case of real multiplication considered above, it is enough to
construct an ample sequence $(\FF_n, n\in\Z)$ of objects in $\CC^{\th}$,
however, not necessarily of the form $\FF_n=F^n(\FF_0)$ for some autoequivalence
$F$. Then the main theorem of \cite{P-coh} will give 
an equivalence $\CC^{\th}\simeq\cohproj A$, where $A$ is the
$\Z$-algebra associated with the sequence $(\FF_n)$, so that
$A_{i,j}=\Hom(\FF_i,\FF_j)$.
The construction of the following theorem provides plenty of
ample sequences in $\CC^{\th}$.

\begin{thm}\label{projthm2} For every $\th\in\R$ there exists
an ample sequence $(\FF_n, n\in\Z)$ in $\CC^{\th}$ such that
all the objects $\FF_n$ are stable.
\end{thm}

\Pf . Clearly, it suffices to consider the case when $\th$ is irrational.
Recall that all vectors $v_{\FF}$ for $\FF\in\CC^{\th}$ belong to
the half-plane $H=H_{\th}=\{(x_1,x_2)|\ x_1-\th x_2>0\}\sub \R^2$. 
Moreover, for every primitive vector $v\in H\cap\Z^2$ there exists a stable
object $\FF\in\CC^{\th}$ with $v_{\FF}=v$. Now let us choose
a sequence of primitive vectors $v_n=(d_n,r_n)\in H\cap\Z^2$
such that $r_n>0$ for $n<<0$ and $\lim_{n\to-\infty}\mu_n=\th$,
where $\mu_n=d_n/r_n$. In other words, we want the ray $\R_{\ge 0}v_n$
to approach $\R_{\ge 0}(\th,1)$ as $n\to-\infty$. 
Note that since $d_n-\th r_n>0$ we necessarily have $\mu_n>\th$ for $n<<0$.
In addition we can make this choice
in such a way that for all $n<<0$ one has $\mu_n-\th\ge r_n^{-1}$.
Indeed, we can first choose $r_n$ for $n<<0$ to be a sequence
of prime numbers such that $\lim_{n\to-\infty}r_n=+\infty$.
Then after picking any sequence $d_n$ such that 
$\lim_{n\to-\infty}d_n/r_n=\th$ we can change $d_n$ by $d_n+1$ if necessary
to make $\mu_n-\th\ge r_n^{-1}$ for $n<<0$. Since
$\th$ is not an integer and $d_n/r_n$ tends to $\th$ as $n\to-\infty$,
such a change will leave $d_n$ prime to $r_n$.

Now we claim that if $(\FF_n, n\in\Z)$ is
any sequence of stable objects in
$\CC^{\th}$ with $v_{\FF_n}=v_n$ then conditions (i) and (ii) from
the proof of Theorem \ref{amplethm} are satisfied for every stable
$\FF\in\CC^{\th}$, and therefore the sequence $(\FF_n)$ is ample.
Indeed, condition (i) follows from Lemma \ref{nzhomlem}
since for every $v\in H$ one has $\chi(v_n,v)>0$ for $n<<0$.
Arguing in the same way as in the proof of Theorem \ref{amplethm}
we conclude that it is enough to prove that for every $v\in H\cap\Z^2$ one has 
$\chi(v_n,v)v_n-v\in H$ for $n<<0$.
Let $v=(d,r)$. Then we have to show that
$$(dr_n-d_nr)d_n-d-\th[(dr_n-d_nr)r_n-r]>0$$
for $n<<0$. Assume first that $r\neq 0$ and set $\mu=d/r$.
Then the above inequality can be rewritten as
$$rr_n^2(\mu-\mu_n)(\mu_n-\th)>r(\mu-\th).$$
Note that $r(\mu-\mu_n)=\chi(v_n,v)/r_n>0$ for $n<<0$. 
Hence, our inequality for $n<<0$ is
equivalent to
$$r_n^2(\mu_n-\th)>\frac{\mu-\th}{\mu-\mu_n}.$$
But this follows from the condition that $\mu_n-\th\ge r_n^{-1}$
since $r_n\to+\infty$ as $n\to-\infty$. Similar argument works
in the case $r=0$.
\ed

\end{document}